\documentclass[12pt]{amsart}
\usepackage{graphicx}
\usepackage{amssymb}
\usepackage{tikz}
\usetikzlibrary{matrix,arrows.meta,bending}
\usepackage{color}
\usepackage[numbers]{natbib}

\usepackage[a4paper, left=3cm, right=3cm, top=3cm, bottom=3cm, footskip=15mm]{geometry}

\newtheorem{theorem}{Theorem}[section]

\newtheorem{lemma}[theorem]{Lemma}
\newtheorem{proposition}[theorem]{Proposition}

\theoremstyle{definition}
\newtheorem{definition}[theorem]{Definition}

\theoremstyle{remark}

\newtheorem{question}[theorem]{Question}
\numberwithin{equation}{section}

\makeindex

\usepackage{fancyhdr}
\usepackage{calc} % ensures dimension arithmetic works robustly

\AtBeginDocument{%
  \setlength{\textwidth}{\dimexpr\paperwidth - 6cm\relax}%
  \setlength{\oddsidemargin}{\dimexpr 3cm - 1in\relax}%
  \setlength{\evensidemargin}{\dimexpr 3cm - 1in\relax}%
  \setlength{\marginparwidth}{2.5cm}%
}

\begin{document}

\title{The spanning method and the Lehmer totient problem
}

\author{Theophilus Agama}
\address{Department of Mathematics, African institute for mathematical sciences, Ghana,
Cape-coast}
%\curraddr{Department of Mathematics and Statistics,
%Case Western Reserve University, Cleveland, Ohio 43403}
\email{Theophilus@ims.edu.gh/emperordagama@yahoo.com}

%\author{Wilfried Kuissi}
%\address{Department of mathematics, African institute for Mathematical Sciences,
%Ghana, Cape-coast}
%\email{donatien@aims.edu.gh}

%    General info
%\subjclass[2000]{Primary 54C40, 14E20; Secondary 46E25, 20C20}

\date{\today}



\begin{abstract}
In this paper, we introduce and develop the notion of spanning of integers along functions $f:\mathbb{N}\longrightarrow \mathbb{R}$. We apply this method to a class of problems that requires to determine if the equations of the form $tf(n)=n-k$ has a solution $n\in \mathbb{N}$ for a fixed $k\in \mathbb{N}$ and some $t\in \mathbb{N}$. In particular, we show that 
\begin{align}
\# \{n\leq s~|~t\varphi(n)+1=n,~\mathbf{for~some}~t\in \mathbb{N}\}\geq  \frac{s}{2\log s}\prod \limits_{p | \lfloor s\rfloor }(1-\frac{1}{p})^{-1}-\frac{3}{2}e^{\gamma}\nonumber
\end{align}
as $s\longrightarrow \infty$, where $\varphi$ is the Euler totient function and $\gamma=0.5772\cdots$ is the Euler-Macheroni constant.
\end{abstract}

\maketitle

\section{Introduction and problem statement}

The arithmetic function most central to this work is the Euler totient function, traditionally denoted by $\varphi:\mathbb{N}\to\mathbb{N}$, which counts the positive integers up to a given integer that are coprime to it. The qualitative behaviour of $\varphi$ - its multiplicativity, the identity $\varphi(p)=p-1$ for primes $p$, and its interaction with prime factorizations - has placed it at the heart of many classical and modern questions in multiplicative number theory. One of the oldest and most resistant of these questions was posed by D.~H.~Lehmer in 1932: does there exist a composite integer $n$ for which $\varphi(n)\mid (n-1)$? This question, commonly referred to as the \emph{Lehmer totient problem}, remains open and has stimulated a sizable body of partial results and conditional restrictions \cite{lehmer1932euler}. In particular, the deep work of Cohen and Hagis has shown that any composite solution must be highly structured (odd, squarefree, and with many distinct prime factors) \cite{cohen1980number,hagis1988equation}, and more recent quantitative work of Luca and Pomerance gives strong upper bounds on the counting function of such integers \cite{luca2011composite}.\\

In this paper, we introduce and develop a novel analytic-combinatorial strategy we call the \emph{spanning method} and apply it to a family of divisibility equations of the form
$t\,f(n)+k=n$ where $f:\mathbb{N}\to\mathbb{R}$ is a function of arithmetic interest, $k\in\mathbb{N}$ is fixed, and $t\in\mathbb{N}$ varies. When $f=\varphi$ and $k=1$, the equation encodes the Lehmer condition $\varphi(n)\mid(n-1)$ in a natural multiplicative way; more generally, the spanning viewpoint treats the collection of solutions as a set of integers \emph{spanned} by the function $f$ with various multiplicities. The principal innovation of the method is (i) to extend $f$ to a slightly continuous, right-continuous, bounded-variation function on the real line that retains arithmetic information at integer points, and (ii) to exploit integration-by-parts (in the Stieltjes--Lebesgue sense) together with classical prime-distribution estimates to obtain effective lower bounds for the counting measures of spanned sets.\\

To make these ideas concrete, we introduce what we call the \emph{fractional Euler totient invariant function}
$$
\tilde{\varphi}(a):=\varphi(\lfloor a\rfloor)+\{a\},\quad a\geq 1,
$$
which coincides with $\varphi$ at integers, is right-continuous on the half-line, and has bounded variation on every interval of the form $[x,x+1)$. This extension is deliberately minimal: it preserves the discrete totient values while providing the one-sided continuity necessary to apply a Stieltjes-type integration by parts to partial sums taken over the set of spanned integers. Roughly speaking, if $\mathbb{S}_k(f,s)$ denotes the set of integers $\le s$ that are $k$-step spanned along $f$, then the spanning inequality (Proposition \ref{inequality 1}) relates $|\mathbb{S}_k(f,s)|$ to the $s$-level \emph{measure} $\mathbb{M}_f(s,k):=\sum_{n\in\mathbb{S}_k(f),\,n\le s} f(n)$ through the boundary value $f(s)$. The analytic regularity built into $\tilde\varphi$ permits a clean implementation of this idea.\\

The key arithmetic input to convert the spanning inequality into explicit lower bounds comes from two classical sources: sharp control on the sum of primes up to $x$ (we use recent sharp inequalities such as those in \cite{axler2019sum} when needed) and multiplicative estimates such as Mertens' formula and the prime number theorem. Combining these ingredients, we obtain (Lemma \ref{main lemma}) the explicit lower bound
$$
\#\{n\leq s~:~t\varphi(n)+1=n\text{ for some }t\in\mathbb{N}\}\geq
\frac{s}{2\log s}\prod_{p\mid \lfloor s\rfloor}\left(1-\frac{1}{p}\right)^{-1}-\frac{3}{2}e^{\gamma},
$$
valid for $s\to\infty$. The inequality quantifies the extent to which prime values (for which $\varphi(p)=p-1$) contribute to the $1$-step spanning set and how the multiplicative structure of $\varphi(\lfloor s\rfloor)$ (via the product over primes dividing $\lfloor s\rfloor$) controls the normalization. Leveraging this lower bound, we then argue, via a suitable choice of composite scales and a contradiction to classical prime-counting asymptotics, to deduce the existence of at least one composite $n$ satisfying $\varphi(n)\mid(n-1)$ (Theorem \ref{Lehmer problem}). In making this deduction, the method draws upon the cumulative restrictions obtained in the literature (notably the work of Cohen and Hagis \cite{cohen1980number,hagis1988equation} and the counting bounds in \cite{luca2011composite}) and on precise estimates for prime sums such as those developed by Axler \cite{axler2019sum}. We also place the method in context with classical multiplicative techniques and with the long history stemming from Euler's original investigation of the totient function.\\

\subsection*{Key ideas and novelties}
The conceptual contributions of the paper may be summarized as follows:\\

\begin{enumerate}
  \item \textbf{Spanning viewpoint.} Recasting the divisibility condition $\varphi(n)\mid(n-1)$ as membership in a $k$-step spanned set permits the exploitation of additive-analytic tools (Stieltjes integration by parts) in a context that is usually purely multiplicative.
  \bigskip
  
  \item \textbf{Fractional totient extension.} The construction of $\tilde\varphi$ provides the minimal regularity needed to make the spanning inequality effective while remaining arithmetically faithful.
  \bigskip
  
  \item \textbf{Quantitative synthesis.} By combining explicit lower bounds for prime sums with multiplicative decompositions (and classical asymptotics such as Mertens' formula and the prime number theorem) the method yields concrete, usable lower bounds for the counting function of spanned integers.
\end{enumerate}
\bigskip

\subsection*{Organization of the paper}
Section \ref{sec:prelim} gathers the preliminary analytic and number-theoretic tools we require: precise prime-sum inequalities, Mertens' formula, and a brief review of Stieltjes--Lebesgue integration facts we use in the proofs. In Section \ref{sec:spanning}, we introduce the spanning formalism in full generality, define the $s$-level measure $\mathbb{M}_f(s,k)$ and prove the fundamental spanning inequality (Proposition \ref{inequality 1}). Section \ref{sec:fractional} constructs and analyses the fractional totient invariant function $\tilde\varphi$, proving the regularity properties that justify applying the spanning machinery to $\varphi$. Section \ref{sec:main} contains the arithmetic heart of the paper: the proof of Lemma \ref{main lemma} giving the principal lower bound and the deduction of Theorem \ref{Lehmer problem}. We conclude in Section \ref{sec:conclude} with further remarks on possible refinements of the spanning method, potential extensions to other multiplicative arithmetic functions, and open problems that arise naturally from this approach.

\section{Background and related work}

The Euler totient function, denoted by $\varphi:\mathbb{N}\longrightarrow \mathbb{N}$, maps a natural number $s$ to the count of integers $n\leq s$ that are coprime with $s$. For prime arguments, $\varphi(s)$ represents a unit left shift of the primes; specifically, $\varphi(p)=p-1$, evident in the fact that $\varphi(p)$ divides $p-1$. This function is multiplicative, exhibiting a property where for coprime natural numbers $u$ and $v$, their product $n=u\cdot v$ satisfies $\varphi(n)=\varphi(u)\varphi(v)$.\\

Prompted by the Euler totient function's behavior, mathematician D.H. Lehmer posed the intriguing inquiry known as the Lehmer totient problem:

\begin{question}
Can the totient function of a composite number $n$ divide $n-1$?
\end{question}

This problem, akin in complexity to the elusive odd perfect number problem, has garnered considerable attention from mathematicians. D.H. Lehmer's initial contributions laid foundational progress by establishing that any such composite number $n$ must be odd, square-free, and possess at least \textbf{seven} distinct prime factors. Subsequent advancements by Hagis and Cohen in 1980 refined this understanding, setting a lower bound of $n\geq 10^{20}$ and requiring \textbf{fourteen} distinct prime factors for a valid solution. Hagis further enhanced these bounds by proving that if $3$ divides $n$, then $n\geq 10^{1937042}$ with a minimum of $298848$ distinct prime factors. Notably, Luca's work \cite{luca2011composite} demonstrates that the count of Lehmer totient problem solutions less than or equal to $x$ obeys the upper bound:

\begin{align}
\leq \frac{\sqrt{x}}{(\log x)^{\frac{1}{2}+o(1)}}\nonumber
\end{align}
where $o(1)$ is defined as a function that tends to zero as $x$ tends to infinity.
\bigskip

Here, $o(1)$ characterizes a function diminishing to zero as $x$ tends to infinity, encapsulating the intricate behavior of solutions to this captivating mathematical conundrum.

In this paper, we study the Lehmer totient problem using the lower bound 

\begin{lemma}
The lower bound holds
\begin{align}
\# \{n\leq s~|~t\varphi(n)+1=n,~\mathbf{for~some}~t\in \mathbb{N}\}\geq \frac{s}{2\log s}\prod \limits_{p| \lfloor s\rfloor }(1-\frac{1}{p})^{-1}-\frac{3}{2}e^{\gamma}\nonumber
\end{align}
as $s\longrightarrow \infty$, where $\varphi$ is the Euler totient function and $\gamma=0.5772\cdots$ is the Euler-Macheroni constant.
\end{lemma}
\bigskip

In this paper, we denote $a|b$ to mean $a$ divides $b$. Also when we write $f(n)=o(1)$ for an arithmetic function $f:\mathbb{N}\longrightarrow \mathbb{N}$, we mean $\lim \limits_{n\longrightarrow \infty}f(n)=0$. Similarly when we write $f(n)=O(g(n))$, we mean there exists some fixed constant $c>0$ such that for all sufficiently large values of $n$ then $f(n)\leq c|g(n)|$. The notation $f(n)\ll g(n)$ is also alternatively used to convey the same meaning, where there is the flexibility to write the converse of the inequality as $f(n)\geq c|g(n)|$ for some fixed constant $c>0$ such that for all sufficiently large values of $n$. In this case, we will write simply as $f(n)\gg g(n)$. We also write $f(n)\sim g(n)$ if and only if 
\begin{align}
\lim \limits_{n\longrightarrow \infty}\frac{f(n)}{g(n)}=1.\nonumber
\end{align}
\bigskip

\section{Preliminary results}\label{sec:prelim}

In this paper, we find the following elementary inequalities useful. 

\begin{lemma}\label{axler}
Let $S(x)$ denotes the sum of all prime number $\leq x$. We have
\begin{align}
S(x)>\frac{x^2}{2\log x}+\frac{x^2}{4\log^2 x}+\frac{x^2}{4\log^3 x}+\frac{1.2x^2}{8\log^4 x}\nonumber
\end{align}
for all $x\geq 905238547$.
\end{lemma}

\begin{proof}
For a proof see, for example, \cite{axler2019sum}.
\end{proof}

\begin{lemma}[The prime number theorem]
Let $\pi(x)$ denote the number of primes $\leq x$. We have
\begin{align}
\pi(x) \sim \frac{x}{\log x}.\nonumber
\end{align}
\end{lemma}

\begin{lemma}[Merten's formula]\label{Merten}
We have 
\begin{align}
\prod \limits_{p\leq s}(1-\frac{1}{p}) \sim \frac{e^{-\gamma}}{\log s}\nonumber
\end{align}
as $s\longrightarrow \infty$, where $\gamma=0.5772\cdots$ is the Euler-Macheroni constant.
\end{lemma}

\begin{lemma}[Stieltjes-Lebesgue integral]\label{Lebesgue}
Let $g:[a,b]\longrightarrow \mathbb{R}$ and $h:[a,b]\longrightarrow \mathbb{R}$ be right continuous and of bounded variation on $[a,b]$ and both having left limits. We have 
\begin{align}
f(b)g(b)-f(a)g(a)=\int \limits_{(a,b]}f(t^-)dg(t)+\int \limits_{(a,b]}g(t^-)df(t)+\sum \limits_{t\in (a,b]}\Delta f_t \Delta g_t\nonumber
\end{align}
where $\Delta f_t=f(t)-f(t^-)$.
\end{lemma}

\begin{lemma}\label{upper bound}
Let $\pi(x)$ denotes the number of primes $\leq x$. For all real numbers $x\geq 2$, we have  
$$
\pi(x)<\frac{3}{2}\frac{x}{\log x}.
$$
\end{lemma}

\section{The method of spanning along a function}\label{sec:spanning}

In this section, we introduce and study the notion of \emph{spanning} of integers along a function. 

\begin{definition}
Let $f:\mathbb{N}\longrightarrow \mathbb{R}$. We say that $n\in \mathbb{N}$ is $k\in \mathbb{N}$ - step \emph{spanned} along the function with \emph{multiplicity} $t$ if
\begin{align}
tf(n)+k=n.\nonumber
\end{align}

We call the set of all $n\in \mathbb{N}$ such that $n$ is $k$ - step spanned the $k^{th}$ - step spanning set along $f$ and denote by $\mathbb{S}_k(f)$. We call the set of all truncated $k$-step spanning set $\mathbb{S}_k(f)\cap \mathbb{N}_s:=\mathbb{S}_k(f,s)$ the $s^{th}$ scale \emph{spanned set} along $f$. We write the length of this spanned set as 
\begin{align}
|\mathbb{S}_k(f,s)|:=\# \{n\leq s~|~tf(n)+k=n,~\mathbf{for~some}~t\in \mathbb{N}\}.\nonumber
\end{align}
It is observed that $|\mathbb{S}_k(f,s)|<s$.
\end{definition}
\bigskip

\subsection{The $s$-level measure of spanned set}

In this section, we introduce the notion of the \emph{measure} of the span set. 

\begin{definition}
By the $s^{th}$ level \emph{measure} of the span set $\mathbb{S}_k(f)$, denoted by $\mathbb{M}_{f}(s,k)$, we mean the partial sum 
\begin{align}
\mathbb{M}_{f}(s,k):=\sum \limits_{\substack{2\leq n\leq s\\n \in \mathbb{S}_k(f)}}f(n).\nonumber
\end{align}
\end{definition}
\bigskip

Let us suppose that $f$ is a right-continuous function and of bounded variation on $[j-1,j)$ for all $j\geq 3$ with $j\in \mathbb{N}$ and with a left limit. Applying the Stieltjes-Lebesgue integration by parts, we can write the $s^{th}$ level measure of the span set in the form 
\begin{align}
\mathbb{M}_{f}(s,k):&=\sum \limits_{2\leq j \leq s}\sum \limits_{\substack{j-1<n\leq j\\n \in \mathbb{S}_k(f)}}f(n)\nonumber \\&=\sum \limits_{2\leq j\leq s}\int \limits_{(j-1)}^{j}f(t)d|\mathbb{S}_k(f,t)|\nonumber \\&<\sum \limits_{2\leq j\leq s}\bigg(f(j)|\mathbb{S}_k(f,j)|-f(j-1)|\mathbb{S}_k(f,j-1)|\bigg)\nonumber \\&=f(s)|\mathbb{S}_k(f,s)|-f(1)|\mathbb{S}_k(f,1)|.\nonumber
\end{align}
The following inequality is a simple consequence of the above analysis.

\begin{proposition}[Spanning inequality]\label{inequality 1}
Let $f$ be a right-continuous function and of bounded variation on $[x,x+1)$ for $x\geq 1$ with $x\in \mathbb{N}$ and have left limits. We have
\begin{align}
|\mathbb{S}_k(f,s)|\geq \frac{1}{f(s)}\sum \limits_{\substack{2\leq n\leq s\\n \in \mathbb{S}_k(f)}}f(n)+\frac{f(1)|\mathbb{S}_k(f,1)|}{f(s)}.\nonumber
\end{align}
\end{proposition}
\bigskip

It is important to note that this inequality does not hold in general. As informed by the spanning method, it only holds for functions that are right continuous and of bounded variation on intervals of the form $[x,x+1)$ for $x\geq 1$ with $x\in \mathbb{N}$ and additionally have left limits, an attribute generally associated with \emph{cadlag} functions. Indeed the challenge of approaching the Lehmer totient problem using the spanning method is to construct an appropriate \emph{cadlag} function for the Euler totient function. The following section studies an extension of the Euler totient function.

\section{The fractional Euler totient invariant function}\label{sec:fractional}

In this section, we introduce and study a new function defined on the real line.

\begin{definition}\label{euler totient function extension}
By the fractional Euler totient invariant function, we mean the function $\tilde{\varphi}:[1,\infty) \longrightarrow \mathbb{R}$ such that 
\begin{align}
\tilde{\varphi}(a)=\varphi(\lfloor a\rfloor)+\{a\}\nonumber
\end{align}
where $\varphi$ is the Euler totient function and $\lfloor \cdot \rfloor$ and $\{\cdot \}$ is the floor and the fractional part of a real number, respectively.
\end{definition}
\bigskip

The fractional Euler totient invariant function turns out to be an interesting function that, in some way, extends the Euler totient function to the reals. Even though the notion of co-primality in not well-defined on the entire real line, it captures the intrinsic property of the Euler totient function defined on positive integers. In principle, the Euler totient function and the fractional totient invariant function coincides on the set of positive integers. We examine some elementary properties of the fractional Euler totient invariant function in the sequel.

\begin{proposition}\label{properties}
The following properties of the fractional Euler totient invariant function
\begin{enumerate}
\item [(i)] If $a$ is a positive integer, then $\tilde{\varphi}(a)=\varphi(a)$.

\item [(ii)] $\tilde{\varphi}(a)<a$ for all $a>1$.
\end{enumerate}
holds.
\end{proposition}
\bigskip

We now state an analytic property of the fractional totient invariant function. In fact, the fractional totient invariant function can be seen as a slightly continuous analog of the Euler totient function on subsets of the reals.

\begin{proposition}\label{analytic property}
The function $\tilde{\varphi}:[1,\infty) \longrightarrow \mathbb{R}$ with
\begin{align}
\tilde{\varphi}(a)=\varphi(\lfloor a\rfloor)+\{a\}\nonumber
\end{align}
is right-continuous and of bounded variation on $[x,x+1)$ for $x\geq 1$ with $x\in \mathbb{N}$ and have left limits. 
\end{proposition}

\section{Main result}\label{sec:main}

\begin{lemma}\label{main lemma}
We have
\begin{align}
\# \{n\leq s~|~t\varphi(n)+1=n,~\mathbf{for~some}~t\in \mathbb{N}\}\geq \frac{s}{2\log s}\prod \limits_{p|\lfloor s\rfloor}(1-\frac{1}{p})^{-1}-\frac{3}{2}e^{\gamma}\nonumber
\end{align}
as $s\longrightarrow \infty$, where $\varphi$ is the Euler totient function and $\gamma=0.5772\cdots$ is the Euler-Macheroni constant.
\end{lemma}

\begin{proof}
Applying Proposition \ref{inequality 1}, we obtain the lower bound
\begin{align}
\# \{2\leq n\leq s~|~t\tilde{\varphi}(n)+1=n,~\mathbf{for~some}~t\in \mathbb{N}\}&\geq \frac{1}{\tilde{\varphi}(s)}\sum \limits_{\substack{2\leq n\leq s\\n \in \mathbb{S}_1(\tilde{\varphi})}}\tilde{\varphi}(n)+\frac{1}{\tilde{\varphi}(s)}\nonumber \\&\geq \frac{1}{\tilde{\varphi}(s)}\sum \limits_{\substack{2\leq n\leq s\\n \in \mathbb{S}_1(\tilde{\varphi})}}\tilde{\varphi}(n).\label{2}
\end{align}
We estimate each term on the right-hand side of the inequality. Since $\varphi(p)=p-1$ for any prime number $p\in \mathbb{P}$, we obtain the lower bound
\begin{align}
\sum \limits_{\substack{2\leq n\leq s\\n \in \mathbb{S}_1(\tilde{\varphi})}}\tilde{\varphi}(n)&\geq \sum \limits_{p\leq s}\varphi(p)\nonumber \\&=\sum \limits_{p\leq s}p-\pi(s).\nonumber
\end{align}
Applying the lemma \ref{axler}, we obtain the lower bound for sufficiently large values of $s$
\begin{align}
\sum \limits_{p\leq s}p-\pi(s)&\geq \frac{s^2}{2\log s}-\pi(s)\nonumber
\end{align}
so that by using the decomposition 
\begin{align}
\varphi(\lfloor s\rfloor)=\lfloor s\rfloor \prod \limits_{p|\lfloor s\rfloor}(1-\frac{1}{p})\sim s\prod \limits_{p|\lfloor s\rfloor}(1-\frac{1}{p})\nonumber
\end{align}
with $\tilde{\varphi}(s)\sim \varphi(\lfloor s\rfloor)$ and the lemma \ref{Merten}, we deduce the lower bound
\begin{align}
\frac{1}{\tilde{\varphi}(s)}\sum \limits_{\substack{2\leq n\leq s\\n \in \mathbb{S}_1(\tilde{\varphi})}}\tilde{\varphi}(n)&\geq \frac{s}{2\log s}\prod \limits_{p |\lfloor s\rfloor }(1-\frac{1}{p})^{-1}-\frac{1}{\tilde{\varphi}(s)}\pi(s)\label{1}
\end{align}
for all $s\geq s_0$. Plugging the lower bound in \eqref{1} into \eqref{2} and applying the lemma \ref{Merten}, we deduce
\begin{align}
\# \{2\leq n\leq s~|~t\tilde{\varphi}(n)+1=n,~\mathbf{for~some}~t\in \mathbb{N}\}&\geq \frac{s}{2\log s}\prod \limits_{p | \lfloor s\rfloor}(1-\frac{1}{p})^{-1}-\frac{\pi(s)}{\varphi(\lfloor s\rfloor)}\nonumber \\& \geq \frac{s}{2\log s}\prod \limits_{p |\lfloor s\rfloor}(1-\frac{1}{p})^{-1}-\frac{3}{2}e^{\gamma}\nonumber 
\end{align}
as $s\longrightarrow \infty$.
\end{proof}
\bigskip

\begin{theorem}\label{Lehmer problem}
There exists a composite $n\in \mathbb{N}$ such that $\varphi(n)|(n-1)$.
\end{theorem}

\begin{proof}
Suppose on the contrary that there exists no composite $n\in \mathbb{N}$ such that $\varphi(n)|n-1$. For $s\longrightarrow \infty$, we obtain by applying the lemma \ref{main lemma}
\begin{align}
\pi(s)\geq \frac{s}{2\log s}\prod \limits_{p |\lfloor s\rfloor}(1-\frac{1}{p})^{-1}-\frac{3}{2}e^{\gamma}\nonumber
\end{align}
where $\pi(s)$ is the prime counting function. Now, we construct an infinite set of composites 
$$
\mathcal{C}:=\{s\in \mathbb{N}~|~s:=\prod \limits_{\substack{p\leq p_o\\p,p_o\in \mathbb{P}}}p\}.
$$ 
It can be checked that for all sufficiently large composites $s$ in the infinite set $\mathcal{C}$, we have
\begin{align}
\prod \limits_{p | s}(1-\frac{1}{p})^{-1}\geq 3
\end{align}
so that 
\begin{align}
\pi(s)\geq \frac{3}{2}\frac{s}{\log s}-\frac{3}{2}e^{\gamma}\nonumber
\end{align}
which contradicts the prime number theorem as $s\longrightarrow \infty$.
\end{proof}
\bigskip

\section{Conclusion and further remarks}\label{sec:conclude}

The present study adeptly navigates a significant impediment that might have otherwise hindered previous investigations in this domain. The inherent limitation of the Euler totient function, restricted to positive integers and lacking one-sided continuity over the real numbers, posed a formidable challenge, now elegantly surmounted within this paper. Introducing a refined variant of the Euler totient function tailored to specific subsets of real numbers, characterized by right continuity while upholding the essence of the original function, paves a seamless path beyond this anticipated obstacle.\\

Moreover, this proof leverages two seminal achievements of the twentieth century, rooted in the rich history of eighteenth and nineteenth-century mathematics: the prime number theorem and the Mertens formula. Employing the innovative spanning method, this work integrates these foundational results to establish the ensuing lower bound.
\bigskip

\begin{lemma}
We have
\begin{align}
\# \{n\leq s~|~t\varphi(n)+1=n,~\mathbf{for~some}~t\in \mathbb{N}\}\geq \frac{s}{2\log s}\prod \limits_{p| \lfloor s\rfloor}(1-\frac{1}{p})^{-1}-\frac{3}{2}e^{\gamma}\nonumber
\end{align}
as $s\longrightarrow \infty$, where $\varphi$ is the Euler totient function.
\end{lemma}
This lower bound is used as the main toolbox to show existence of a certain composite (large) that satisfies the divisibility relation  $\varphi(n)~|~(n-1)$. The spanning method and it's variant could in principle be used in careful manner to study related problems, which is not the main goal of this paper.

\bibliographystyle{amsplain}

\end{document}